\documentclass[12pt]{article}
\usepackage{geometry,amsmath,amssymb,graphicx,bbm,xypic,relsize,exscale,amscd,amsthm,url}
\newtheorem{theorem}{Theorem}

\newtheorem{conjecture}[theorem]{Conjecture}

\newcommand{\tmtextit}[1]{{\itshape{#1}}}

\begin{document}

\title{Metaplectic Theta Functions and Global Integrals}
\author{Solomon Friedberg and David Ginzburg
\footnote{This work was supported by the US-Israel Binational Science Foundation,
grant number 201219, and by the National Security Agency, grant number
H98230-13-1-0246 (Friedberg).}\\
\small Department of Mathematics, Boston College, Chestnut Hill MA 02467-3806, USA\\
\small School of Mathematical Sciences, Tel Aviv University, Ramat Aviv, Tel Aviv 6997801,
Israel}
\date{}
\maketitle

\centerline{\it To Steve Rallis, in memoriam}

\begin{abstract} We convolve a theta function on an $n$-fold
cover of $GL_3$ with an automorphic form on an $n'$-fold cover of $GL_2$
for suitable $n,n'$.
To do so, we induce the theta function to the $n$-fold cover of $GL_4$ and use a Shalika integral.
We show that in particular when $n=n'=3$  this construction gives a new Eulerian integral
for an automorphic form on the 3-fold cover of $GL_2$ (the first such integral was given
by Bump and Hoffstein), and when $n=4$, $n'=2$, it gives a Dirichlet series
with analytic continuation and functional equation that involves both the Fourier
coefficients of an automorphic form of half-integral weight and quartic
Gauss sums.  The analysis of these cases is based on the uniqueness of the Whittaker
model for the local exceptional representation.
The constructions studied here may be put in the context of a larger family
of global integrals which are constructed using automorphic representations on
covering groups. We  sketch this wider context and some related conjectures.
\end{abstract}

\noindent
MSC2010 Subject Classification:  Primary 11F70; Secondary 11F27, 11F37, 11F55
\smallskip


\section{Introduction}
In this paper we initiate the study of a certain family of global
integrals which are constructed using automorphic representations on
covering groups. Beside their intrinsic interest,
such constructions have the potential to
relate the values of Whittaker coefficients of
the theta representations defined on covering groups of
general linear groups of different ranks, a sort of metaplectic
version of descent.  To begin, we explain the general
context for these constructions.

Let $F$ be a global field, $\mathbb{A}$ be the adeles of $F$, and
$GL_m^{(n)}(\mathbb{A})$ be an n-fold cover of $GL_m(\mathbb{A})$,
defined when $F$ has a full set of $n$-th roots of unity.  (For convenience
we shall assume that $F$ has a full set of $2n$-th roots of unity below.) There is
more than one such cover, but the 2-cocyles describing different covers are related
by twists, and agree on $SL_m(\mathbb{A})$. Let
$\Theta_m^{(n)}$ denote the theta representation on
$GL_m^{(n)}(\mathbb{A})$, constructed by Kazhdan and Patterson
\cite{K-P} via residues of Eisenstein series. For $r\geq2$, let
$E_{\Theta_{2r-1}^{(n)}}(g,s,f_s)$ denote the Eisenstein series
defined on $GL_{2r}^{(n)}(\mathbb{A})$ obtained by parabolically
inducing $\Theta_{2r-1}^{(n)}$.  Since the representation
$\Theta_{2r-1}^{(n)}$ is associated with the unipotent orbit
indexed by the partition $((r+1)(r-2))$, this Eisenstein series is associated with the
unipotent orbit $((r+2)(r-2))$.  

Let $U_{2r}$ denote the unipotent radical of the parabolic subgroup
of $GL_{2r}$ whose Levi part is $GL_2^r$, and $\psi_{U_{2r}}$ denote
a character of $U_{2r}(F)\backslash U_{2r}(\mathbb{A})$ whose
stabilizer inside $GL_2^r$ is the group $GL_2$ embedded diagonally.
The group $U_{2r}(\mathbb{A})$ embeds canonically in $GL_{2r}^{(n)}(\mathbb{A})$
by the trivial section $\mathbf{s}(h)=(h,1)$.
Then the Fourier coefficient
\begin{equation*}
\int\limits_{U_{2r}(F)\backslash U_{2r}(\mathbb{
A})}E_{\Theta_{2r-1}^{(n)}}(\mathbf{s}(u)g,s,f_s)\,\psi_{U_{2r}}(u)\,du\notag
\end{equation*}
is attached to the unipotent orbit $(r^2)$. 
Let $Z^n(\mathbb{A})$ denote the group of scalar matrices $\lambda I_{2}$ in $GL_{2}(\mathbb{A})$ 
such that $\lambda\in\mathbb{A}^\times$ is an $n$-th power.  Then $\mathbf{s}(Z^n(\mathbb{A}))$ is a central subgroup of  
$GL_{2}^{(n)}(\mathbb{A})$.
Let $\pi$ denote a
cuspidal representation of $GL_2^{(k)}(\mathbb{A})$ for certain $k$ with trivial central character. 
Then we form the global integral
\begin{multline}\label{four2}
I(\varphi_\pi,s,f_s)=\\
\int\limits_{Z^n(\mathbb{A})\,GL_2(F)\backslash GL_2(\mathbb{A})}
\int\limits_{U_{2r}(F)\backslash U_{2r}(\mathbb{
A})}\varphi_\pi(\mathbf{s}(g))\,E_{\Theta_{2r-1}^{(n)}}(\mathbf{s}(ug),s,f_s)\,\psi_{U_{2r}}(u)\,du
\end{multline}
where $\varphi_\pi$ is a vector in the space of $\pi$. Here $k$ and the choice of cocycles up to twisting are
chosen so that the cover splits in the integrand on the embedded $GL_2$.  (Even though we mod out by $Z^n(\mathbb{A})$
instead of $Z^1(\mathbb{A})$, using the strong approximation theorem it is not difficult to see that the integral converges.)

Experience with integrals involving Eisenstein series and automorphic
forms suggests that they are most
likely to be Eulerian when the dimension equation is satisfied (see for example Ginzburg \cite{Gi}.)
This is the case here.  In other words,
for this global integral, we have the identity
$$\text{dim}\ GL_2 - \text{dim}\ Z^n+ \text{dim}\ U_{2r}=\text{dim}\ \pi+\text{dim}\
E_{\Theta_{2r-1}^{(n)}},$$ 
where the dimensions on the right-hand side denote the Gelfand-Kirillov dimension as in \cite{Gi}.
Indeed, from the remarks on unipotent orbits
above this identity is equivalent to
$$3+\frac{1}{2}\text{dim}\ (r^2)=1+\frac{1}{2}\text{dim}\ ((r+2)(r-2)),$$ which
is easily verified using, for example \cite{Gi}, Section 2.  

In this paper we begin the study of the following

\begin{conjecture}\label{conje1}
The integral $I(\varphi_\pi,s,f_s)$ has a meromorphic continuation
to the full complex plane, and:\\
{\bf 1)}\ If $n<r+1$, the integral is identically zero for all
$s$.\\
{\bf 2)}\ If $n=r+1$, the integral is Eulerian, and represents the
partial degree two $L$-function $L^S(\tau(\pi),s)$, where
$\tau(\pi)$ is the lift of $\pi$ to $GL_2(\mathbb{A})$.\\
{\bf 3)}\ If $n>r+1$ the integral is not Eulerian and represents a
certain Dirichlet series which, in the domain $\text{Re} (s)>1/2$,
can have at most a simple pole at $s=\frac{n+1}{2n}$.
\end{conjecture}

The second part of this conjecture may be regarded as an
extension of the philosophy of Bump and Hoffstein \cite{B-H0, B-H, B-H2},
who proposed that such an $L$-function could be obtained by ``convolving"  $\varphi_\pi$ with
$\Theta_n^{(n)}$, generalizing the construction of Shimura \cite{S} for $n=2$.  Conjecture~\ref{conje1}
indicates a way to construct such $L$-functions by using theta functions on higher rank groups.
Here a ``convolution" is an adaption of a Rankin-Selberg
integral to the covering group. We note that there is often more than one integral 
representing a given Langlands $L$-function,
and it is not apparent that all extensions of such integrals to the metaplectic group will give the same
results.  Also, since $\varphi_\pi$ does not typically have a {\sl unique} Whittaker model, any such
integral is an example of the class of integrals whose study was initiated by Piatetski-Shapiro and Rallis
in their seminal paper
\cite{PS-R}.  See also Bump and Friedberg \cite{B-F2} who connect the approach of \cite{PS-R} to
the Bump-Hoffstein Conjecture, and Suzuki \cite{Su} for progress on the
Bump-Hoffstein conjecture.

The main result of this paper is an analysis of the global integral (\ref{four2}) in
the case when $r=2$ and $n\leq4$.  We obtain a Dirichlet series for
covers of all degrees but not enough is known about the Whittaker coefficients of higher
theta functions to compute the series for $n>4$.  However, for $n=3$ we show that the
integral is Eulerian and express the local factors in terms of the Hecke eigenvalues
for $\varphi_\pi$.
A different Eulerian integral which is a convolution of $\varphi_\pi$ with $\Theta_3^{(3)}$
was first given by Bump and Hoffstein \cite{B-H2}.
It is interesting that the integral and resulting Dirichlet series arising from the convolution presented here are
fundamentally different from the ones given in \cite{B-H2}, but ultimately represent the same $L$-function.

Let us mention two possible applications of this construction. The first is
in the case $n=r+1$. In this case, Conjecture \ref{conje1} implies that the
integral \eqref{four2} is a holomorphic function in $s$. The
Eisenstein series has a simple pole at $s=\tfrac{r+2}{2(r+1)}$, whose
residue is the representation $\Theta_{2r}^{(r+1)}$. Therefore
replacing the Eisenstein series by
$\Theta_{2r}^{(r+1)}$ in \eqref{four2}, we expect to get zero for all choices of data.
This suggests that the function
\begin{equation}\label{four3}
f(g)=\int\limits_{U_{2r}(F)\backslash U_{2r}(\mathbb{A})}
\theta_{2r}^{(r+1)}(ug)\psi_{U_{2r}}(u)du
\end{equation}
will not have a constituent which is cuspidal. In fact we have
\begin{conjecture}\label{conje2}
The automorphic representation of $GL_2^{(r+1)}(\mathbb{A})$ generated
by all the functions $f(g)$ is $\Theta_2^{(r+1)}$.
\end{conjecture}

Note that if we compare the Whittaker coefficients of both sides
of  \eqref{four3} we obtain a relation between the
Whittaker coefficients of $f$ and the Whittaker coefficients of
other theta functions. For example, when $r=2$, we obtain the
following local identity
\begin{equation*}
q^{n/6}\overline{\tau_{3,2,f}(p^n)}=\tau_{3,3}(1,p^n)+q^{1/6}G_1^{(3)}(p)\tau_{3,3}(1,p^{n-1})
\end{equation*}
Here $\tau_{3,2,f}=W_f\delta_B^{-1/2}$ denotes the normalized
Whittaker coefficient of $f$ at the place $p$, and  the other
notations are as in Hoffstein \cite{Ho}.
It follows from the well known
formulas for the function $W_{\theta_3^{(3)}}$ (see for example
Bump and Hoffstein \cite{B-H0}, \cite{B-H})
that $W_f=W_{\theta_2^{(3)}}$.
More broadly, Conjecture \ref{conje2} implies that the
mapping given by \eqref{four3} is  a descent map in the
sense of Ginzburg, Rallis and Soudry \cite{GRS}.  Such descent
constructions have not been given for general metaplectic covers.

A second situation of interest is $n=r+2$. 
In this case the integral should reduce to one involving
a theta representation that has a unique Whittaker model, and
whose Whittaker coefficients may be computed in terms of {$n$-th}
order Gauss sums. These integrals then give Dirichlet series with
continuation whose coefficients involve both the Fourier
coefficients of $\varphi_\pi$ and arithmetic data.  We illustrate
that here when $r=2$, $n=4$, and where $\pi$ is an automorphic
representation on the double cover, that is, one corresponding to an
automorphic form of half-integral weight. The resulting series
involves the Fourier coefficients of $\varphi_\pi$, whose squares
are related to the central values of twisted $L$-functions, and quartic
Gauss sums, and possesses analytic continuation.

The rest of this paper is organized as follows.  In Section~\ref{notation} we set out the notation. 
Then in Section~\ref{integral} we introduce the basic integral that we study here: the 
integral~(\ref{four2}) when $r=2$.  This is an integral over the Shalika subgroup of $GL_4$.
In Section~\ref{unfolding}, the integral is unfolded and expressed in terms of the Whittaker coefficients of $\varphi_\pi$ and $\Theta_3^{(n)}$.  This series is then analyzed in Section~\ref{analysis}
when $n=3,4$.  The existence of an Euler product when $n=3$ is established in Theorem~\ref{euler-product}.

\section{Notation}\label{notation}

Fix $n\geq2$.  We work over a number field $F$ containing a full set of $2n$-th roots of unity.
Let $S$ be a finite set of places containing the archimedean ones, the ramified ones,
and enough others such that the ring of $S$-integers $\mathcal{O}_S$ has class number 1.
Let $F_S=\prod_{v\in S} F_v$ and embed $F$ in $F_S$ diagonally.  Then for $r\geq1$
$GL_r(\mathcal{O}_S)$ is a discrete subgroup of $GL_r(F_S)$.

The metaplectic covers  of $GL_r(F_S)$ were constructed by Matsumoto following the
work of Kubota when $r=2$.  Convenient references are Bump-Hoffstein \cite{B-H2}
and Kazhdan-Patterson \cite{K-P}.  Recall that the basic cover is constructed from
embedding $GL_r(F_S)$ into $SL_{r+1}(F_S)$ via $g\mapsto \left(\begin{smallmatrix}g&\\&\det(g)^{-1}
\end{smallmatrix}\right)$ and
restricting the $n$-fold cover of $SL_{r+1}(F_S)$.  The cocycle giving this cover may
then be twisted to obtain covers $\widetilde{GL}_r^{(c)}(F_S)$ where $c\in\mathbb{Z}/n\mathbb{Z}$
(see \cite{K-P}, pg.\ 41); these groups are distinct but all contain the $n$-fold cover of $SL_r(F_S)$
constructed by Matsumoto.
For convenience we write $\tilde{G}_j$ for the group $\widetilde{GL}_j^{(c)}(F_S)$, or $\tilde{G}^{(c)}_j$ if it is important
to identify which of the covers we are using.  The group $SL_r(\mathcal{O}_S)$ embeds in $\tilde{G}_r$
by the map $\iota(\gamma)=(\gamma,\kappa(\gamma))$, where $\kappa$ is the Kubota homomorphism
(see for example Brubaker, Bump and Friedberg \cite{BBF}, Section 4).

The reference \cite{K-P} is adelic while \cite{BBF} uses the above notation.
The two approaches are easily connected.  Indeed,
for each place $v$ of $F$ let $\mathcal{O}_v$ denote the ring of
integers in $F_v$.  Then for primes outside $S$,
the hyperspecial maximal
compact subgroup $GL_r(\mathcal{O}_v)$ embeds canonically in the metaplectic group
(see for example \cite{K-P}, Prop.\ 0.1.2).
Thus the strong approximation
theorem shows that working over $F_S$ is equivalent to working adelically
and unramified outside $S$.

The {\sl theta representation} of concern here, $\Theta_3^{(n)}$ or simply $\Theta$,
is a genuine representation on $\tilde{G}_3^{(c)}$.  It may be constructed globally from residues of Eisenstein series
on $\tilde{G}_3^{(c)}$,
and its local constituents $\Theta_v$ at almost all places
may be obtained as the image of a certain unramified principal series representation under an intertwining operator.
See Kazhdan-Patterson \cite{K-P} for details.
By the results of \cite{K-P}, the representations $\Theta_v$ have
a unique Whittaker model for $n=3$ or $n=4$ (for $n=4$ this requires that $c$ be odd; see \cite{K-P} Corollary
I.3.6).  However, for $n>4$ the model is no longer unique.  By contrast, when $n=2$ all Whittaker coefficients are zero.

\section{The Integral}\label{integral}

Let $\varphi$ be a genuine cuspidal automorphic form on the $n'$-fold cover $\widetilde{GL}_2^{(c')}(F_S)$
of $GL_2(F_S)$ where
\begin{equation}\label{covers-match}
2/n+1/n'\in{\mathbb Z},\qquad 4c+c'\equiv 0 \bmod n.
\end{equation}
We suppose that the finite set $S$ is sufficiently large that $\Theta$ and $\varphi$
are unramified outside $S$.
We now present an integral of $\varphi$
against an Eisenstein series constructed
from $\Theta$.  Though we describe the construction for general covers, our main focus will be the cases of covers of degrees
$n=n'=3$ and $n=4$, $n'=2$, as in those cases the Whittaker coefficients of $\Theta$ are uniquely identified.  Using
this identification, we are then able to completely describe the resulting Dirichlet series, which has
continuation and functional equation as it arises from an integral of an Eisenstein series.

Let $P$ be the standard parabolic subgroup of $GL_4$ of type (3,1); the Levi of $P$ is
isomorphic to $GL_3\times GL_1$.  Let $E_{\Theta}(g,s,f_s)$ denote the Eisenstein series on $\tilde{G}_4$
induced from $\Theta$ using a suitable section $f_s$.  This series is the double residue of the minimal parabolic
Eisenstein series on $\tilde{G}_4$ introduced by Kazhdan and Patterson \cite{K-P}, pg.\ 109; see also
Brubaker, Bump and Friedberg \cite{BBF} where the minimal parabolic
Eisenstein series is constructed using similar notation to that
given here (and where its Whittaker expansion is computed).   
The series $E_\Theta(g,s,f_s)$ may be written as a sum
$$E_\Theta(g,s,f_s)=\sum_{\gamma\in P(\mathcal{O}_S)\backslash SL_4(\mathcal{O}_S)}
f_s(\iota(\gamma) g).\qquad g\in\tilde{G}_4$$
(for example, by taking residues in \cite{BBF}, Eqn.\ (23)).
Here we are taking the central character of $\varphi$ to be trivial but one could relax this assumption by incorporating
a suitable character into the Eisenstein series.

Let $M$ denote the algebraic group of $2\times 2$ matrices
and let $R$ denote the Shalika subgroup of $GL_4$
$$R=\left\{\begin{pmatrix}I_2&M\\&I_2\end{pmatrix}\begin{pmatrix}g&\\&g\end{pmatrix}\mid m\in M, g\in GL_2\right\}.$$
Fix an additive character $\psi$ of $F_S$ of conductor $\mathcal{O}_S$.  For $h\in\mathcal{O}_S$,
the map $\psi_h:R(F_S)\to\mathbb{C}^\times$ given by
$$\psi_h\left(\begin{pmatrix}I_2&M\\&I_2\end{pmatrix}\begin{pmatrix}g&\\&g\end{pmatrix}\right)=\psi\left(h\,\rm{tr}(M)\right)$$
is a character of $R(F_S)$ that is trivial on $R(\mathcal{O}_S)$.

Recall that the standard maximal compact subgroup of $GL_2(F_S)$ is $K=\prod_{v\in S} K_v$
where $K_v=GL_2(\mathcal{O}_v)$ if $v$ is nonarchimedean and $K_v=U_2(\mathbb{C})$ if $v$ is
archimedean (hence necessarily complex).  For $GL_r$ with $r>1$, let $\mathbf{s}:GL_r(F_S)\to \tilde{G}_r$
denote the trivial section $\mathbf{s}(g)=(g,1)$.
We suppose that $\varphi$ and $E_\Theta$ have compatible $K$-types:  the function
$$g\mapsto \varphi(g_1\,\mathbf{s}(g))\,E_\Theta\left(g_2\,\mathbf{s}\begin{pmatrix}g&\\&g\end{pmatrix}\right)$$
is right $K$-invariant for any $g_1$, $g_2$ in the corresponding metaplectic groups.

Let $h\in \mathcal{O}_S$.
Then the following integral is well defined:
\begin{multline*}
I(\varphi,s,f_s)=\\
\int_{Z^n(F_S)R(\mathcal{O}_S)\backslash R(F_S)}\varphi(\mathbf{s}(g))E_\Theta\left(\mathbf{s}\left(\begin{pmatrix}I_2&M\\&I_2\end{pmatrix}\begin{pmatrix}g&\\&g\end{pmatrix}\right)\right)
\psi\left(h\,\rm{tr}(M)\right)\,dM\,dg.
\end{multline*}
(We suppress the dependence of this integral on $h$ from the notation.)
Indeed, the covers match by (\ref{covers-match}) so the integrand
 is $R(\mathcal{O}_S)$-invariant.
This is our main object of study.  

\section{Unfolding}\label{unfolding}

To carry out the unfolding, we suppose that $\text{Re}(s)$ is sufficiently large.
Here and below, unless otherwise specified, all elements of $GL_4(F_S)$ are embedded
in $\tilde{G}_4$ via the trivial section $\mathbf{s}$, and we suppress $\mathbf{s}$ from the notation.
We first substitute the Eisenstein series in the form
$$E_\Theta(g,s,f_s)=\sum_{\gamma\in P(\mathcal{O}_S)\backslash SL_4(\mathcal{O}_S)}
\kappa(\gamma) f_s(\gamma w_0 g).$$
Here we are introducing $w_0$, the long element, for convenience.  Then the cosets are parametrized
by their bottom rows $(D_1,D_2,D_3,D_4)$ modulo $\mathcal{O}_S^\times$.
The 4-tuple $(D_1,D_2,D_3,D_4)$  is
relatively prime, and conversely each relatively prime 4-tuple may be completed to an element of $SL_4(\mathcal{O}_S)$.
The cosets $P(\mathcal{O}_S)\gamma$
such that $\gamma w_0$ is in the big Bruhat cell are those with $D_4\neq0$.
As is usual, the cosets not in the big cell contribute zero to the integral and we 
omit them from now on.

On the cosets with $D_4\neq0$, we have a right action of $R(\mathcal{O}_S)$ that allows us to collapse the sum over bottom
rows to a sum over $\gamma$ with bottom row $(D_1,D_2,0,D_4)$ with $D_1$ and $D_2$ modulo $D_4$
(and $D_4$ modulo $\mathcal{O}_S^\times$)
and in exchange to
unfold the integrals.  When we do so, we obtain
\begin{multline}\label{unfold1}
I(\varphi,s,f_s)=\sum_{\gamma} \kappa(\gamma)\
\int_{Z^n(F_S)B_-(\mathcal{O}_S)\backslash GL_2(F_S)}\int_{Y(\mathcal{O}_S)\backslash Y(F_S)}\int_{X(F_S)}\\
\varphi(g)\,f_s\left(\gamma w_0
\begin{pmatrix}I_2&M\\&I_2\end{pmatrix}\begin{pmatrix}g&\\&g\end{pmatrix}\right)
\psi\left(h\,\rm{tr}(M)\right)\,dM\,dg.
\end{multline}
Here $X$ and $Y$ are each copies of the affine plane
and we have written $M=\left(\begin{smallmatrix}X\\Y\end{smallmatrix}\right)$,
and $B_{-}$ is the lower Borel subgroup of $GL_2$.

The sum over $\gamma$ in (\ref{unfold1}) may be parametrized as a sum over nonzero $D_4\in \mathcal{O}_S$,
$D_4$ modulo $\mathcal{O}_S^\times$, and over $D_1,D_2\in \mathcal{O}_S$, each modulo $D_4$,
and such that $\gcd(D_1,D_2,D_4)=1$.
To go farther, we must find coset representatives and compute the Kubota symbol of such a coset representative.
To do so we follow the approach of Brubaker, Bump and Friedberg \cite{BBF}, generalized to all parabolics in
Brubaker and Friedberg  \cite{BF},
and find coset representatives which are products of embeded $GL_2$'s.
Indeed, we parametrize the cosets by representatives $\gamma=\gamma_1\gamma_2$ with
$$\gamma_1=\begin{pmatrix}1&&&\\&a_1&&b_1\\&&1&\\&c_1&&d_1\end{pmatrix}
\qquad\text{and}\qquad
\gamma_2=\begin{pmatrix}a_2&&&b_2\\&1&&\\&&1&\\c_2&&&d_2
\end{pmatrix}.$$
Note that the bottom row of $\gamma$ is then $(c_2d_1,c_1,0,d_1d_2)$ and the sum
will be over $d_1,d_2,c_1,c_2$ with $d_1, d_2$ modulo $\mathcal{O}_S^\times$ and nonzero,
$c_1$ modulo $d_1d_2$, $\gcd(c_1,d_1)=1$
and $c_2$ modulo $d_2$, $\gcd(c_2,d_2)=1$.
We have $\kappa(\gamma)=\left(\frac{d_1}{c_1}\right)\left(\frac{d_2}{c_2}\right)$ where $(-)$ is the $n$-th power
residue symbol as in \cite{BBF}.

Let $\alpha_1,\alpha_2,\alpha_3$ be
the three simple roots of $GL_4$ corresponding to the standard ordering
(so $\alpha_1(\text{diag}(t_1,t_2,t_3,t_4))=t_1/t_2$, etc.) and for each suitable index set $J$ let
$u_{J}(t)$ be the upper triangular
unipotent matrix with $t$ in position corresponding to the root $\alpha_J:=\sum_{j\in J} \alpha_j$,
and $u^-_J(t)$ be the lower triangular unipotent corresponding to $-\alpha_J$.  Let $h_J(t)$
denote the corresponding diagonal matrix in $GL_4$, so $\alpha_J(h_J(t))=t^2$.
Applying the Bruhat decomposition to $\gamma_1$ and $\gamma_2$,
we see that
$$\gamma w_0=u_{23}(b_1/d_1) h_{23}(d_1^{-1})u^-_{23}(c_1/d_1)u_{123}(b_2/d_2)h_{123}(d_2^{-1})u^-_{123}(c_2/d_2)w_0.$$
Embedding each matrix via the trivial section $\mathbf{s}$,
this equality holds for the metaplectic group, provided one adds the factor $(d_1,c_1)_S\,(d_2,c_2)_S$
on the right-hand side.  Here
$(~,~)_S$ is the product of local Hilbert symbols $\prod_{v\in S}(~,~)_v$.  See \cite{BBF}, Eqn.\ (24).
Moreover, by reciprocity for $n$-th power residue symbols, we have
$$(d_1,c_1)_S\,(d_2,c_2)_S\,\kappa(\gamma)=\left(\frac{c_1}{d_1}\right)\left(\frac{c_2}{d_2}\right).$$

Now $u_{23}(b_1/d_1)$ is in the unipotent radical of $P$ and does not affect $f_s$.  And $u^-_{123}(c_2/d_2)w_0
=w_0u_{123}(c_2/d_2)$.   The factor $u_{123}(c_2/d_2)\in R(F_S)$ may be absorbed into the integral by a variable change
that does not introduce a factor.  We multiply (or, more conceptually, apply the Steinberg relations) to see that
$$u^-_{23}(c_1/d_1)u_{123}(b_2/d_2)=u_{123}(b_2/d_2)u_1(-b_2c_1/d_2d_1)u^-_{23}(c_1/d_1).$$
Now the matrix $u^-_{23}(c_1/d_1)$ may be moved rightward and into $R(F_S)$ and absorbed by a variable
change; this introduces a factor of $\psi( h c_1/d_1d_2)$.
Moreover, the matrix $u_{123}(b_2/d_2)$ is in the unipotent radical of $P$ and does not affect $f_s$.
And the matrix $u_1(-b_2c_1/d_2d_1)$ differs from an element of $R$ by an element of the unipotent
radical of $P$.    That is, we obtain the function $f_s$ evaluated as follows:
\begin{multline*}
f_s\Bigg(\begin{pmatrix}d_1^{-1}d_2^{-2}&&&\\&d_1^{-2}d_2^{-1}&&\\&&d_1^{-1}d_2^{-1}&\\&&&1\end{pmatrix}
\begin{pmatrix}1&-b_2c_1/d_1&&\\&1&&\\&&1&-b_2c_1/d_1\\&&&1\end{pmatrix}w_0 \\  
\begin{pmatrix}I_2&M\\&I_2\end{pmatrix}\begin{pmatrix}d_1d_2&&&\\&d_1d_2&&\\&&d_1d_2&\\&&&d_1d_2\end{pmatrix}
\begin{pmatrix} g&\\&g\end{pmatrix}\Bigg).
\end{multline*}
(In these manipulations we use the formulas of \cite{BBF}, Section 3,
and see that no nontrivial Hilbert symbols arise.)

We now move the second matrix in the argument of $f_s$ to the right and change variables in $M,g$.  This gives
\begin{multline*}
\sum_{\substack{0\neq d_1,d_2\in\mathcal{O}_S/\mathcal{O}_S^\times\\ 
c_1 \text{~mod~} d_1d_2,~(c_1,d_1)=1\\ c_2\text{~mod}^\times d_2}}
\left(\frac{c_1}{d_1}\right)\left(\frac{c_2}{d_2}\right)\psi\left(\frac{h c_1}{d_1 d_2}\right)
\int \varphi^{[d_1^{-1}d_2^{-1}]}\left(\begin{pmatrix}1&\\b_2c_1/d_2&1\end{pmatrix}g\right) \\
f_s\left(\begin{pmatrix}d_1^{-1}d_2^{-2}&&&\\&d_1^{-2}d_2^{-1}&&\\&&d_1^{-1}d_2^{-1}&\\&&&1\end{pmatrix}
w_0\begin{pmatrix}I_2&M\\&I_2\end{pmatrix}\begin{pmatrix}g&\\&g\end{pmatrix}\right)
\psi\left( h\,\text{tr}(M) \right)\,dM\,dg.
\end{multline*}
Here and below $a\text{~mod}^\times c$ means $a$ modulo $c$ with $(a,c)=1$,
$b_2$ is a multiplicative inverse of $c_2$ modulo $d_2$, and for $f\in F^\times$,
$$\varphi^{[f]}(g)=\varphi\left(\mathbf{s}\left(\left(\begin{smallmatrix}f&\\&f\end{smallmatrix}\right)\right)g\right).$$

For $\alpha_j$, $j=1,2,3$, a simple root we let $w_j$ be the corresponding simple reflection and let
$w_{j_1j_2}=w_{j_1} w_{j_2}$, etc.  Factor $w_0=w_{21}w_{32}w_{13}$.  Then $w_{21}\in P(\mathcal{O}_S)$
and may be moved leftward and out of the function $f_s$.  (A cocycle computation is necessary
but in fact no nontrivial Hilbert symbol arises.) Similarly $w_{13}\in R(\mathcal{O}_S)$ and may
be moved rightward.
After a variable change sending $g$ to $w_{13}g$ we obtain
\begin{multline}\label{step1}
\sum_{\substack{0\neq d_1,d_2\in\mathcal{O}_S/\mathcal{O}_S^\times\\
c_1 \text{~mod~} d_1d_2,~(c_1,d_1)=1 \\ c_2\text{~mod}^\times d_2}}
\left(\frac{c_1}{d_1}\right)\left(\frac{c_2}{d_2}\right)\psi\left(\frac{h c_1}{d_1 d_2}\right)\\
\int_{Z^n(F_S)B(\mathcal{O}_S)\backslash GL_2(F_S)}\int_{X(F_S)}
\int_{Y(\mathcal{O}_S)\backslash Y(F_S)}
\varphi^{[d_1^{-1}d_2^{-1}]}\left(\begin{pmatrix}1&b_2c_1/d_2\\&1\end{pmatrix}g\right) \\
f_s\left(\begin{pmatrix}d_1^{-1}d_2^{-1}&&&\\&d_1^{-1}d_2^{-2}&&\\&&d_1^{-2}d_2^{-1}&\\&&&1\end{pmatrix}
w_{32}\begin{pmatrix}I_2&M\\&I_2\end{pmatrix}\begin{pmatrix}g&\\&g\end{pmatrix}\right)
\psi\left( h\,\text{tr}(M) \right)\,dM\,dg.
\end{multline}
Here we now have $M=\left(\begin{smallmatrix} Y\\X\end{smallmatrix}\right)$, and $B$
denotes the upper triangular Borel subgroup of $GL_2$.  Note that after the variable
change we have also moved the long
element of $GL_2$ leftward and out of the argument of $\varphi$.

The next step is to expand $\varphi$ in a Fourier expansion.  Recall we have chosen $S$ so that $\varphi$
is $K_v$-fixed for $v\not\in S$.  This expansion is of the form
\begin{equation}\label{F-expansion}
\varphi\left(g\right)=\sum_{0\neq m\in \mathcal{O}_S}\sum_l \frac{b_l(m)}{|m|^{1/2}}\, W_l\left(\left(\begin{smallmatrix}m&\\&1\end{smallmatrix}\right) g \right).
\end{equation}
Here $|m|$ denotes the cardinality of $\mathcal{O}_S/m\mathcal{O}_S$,
$W_l$ runs over a basis for the {\sl finite-dimensional} space of Whittaker functions
(in particular, each function $W_l$ satisfies
$W_l(\left(\begin{smallmatrix}1&x\\&1\end{smallmatrix}\right)g)=\psi(x)W_l(g)$)
and the Fourier coefficients $b_j(m)$ satisfy $b_j(\epsilon m)=b_j(m)$ for all $\epsilon\in\mathcal{O}_S^\times$.
See Bump-Hoffstein \cite{B-H2}, Eqns.\ (3.13), (3.14).  
For later use, we note that if $f\in F^\times$ then the expansion (\ref{F-expansion}) implies that
$$\varphi^{[f]}(g)=
\sum_{0\neq m\in \mathcal{O}_S}\sum_{j,l} X_{lj}^{[f]} \frac{b_l(m)(m,f)_S}{|m|^{1/2}}
W_j\left(\left(\begin{smallmatrix}m&\\&1\end{smallmatrix}\right)g\right),$$
where the complex numbers 
$X_{lj}^{[f]}$ are defined by
$$W_l(\mathbf{s}\left(\left(\begin{smallmatrix}f&\\&f\end{smallmatrix}\right)\right)g)=
\sum_j X_{lj}^{[f]} W_j(g).$$

We now use the Iwasawa decomposition on $GL_2(F_S)$ to replace $g$ in the integral by 
$\left(\begin{smallmatrix}1&x\\&1\end{smallmatrix}\right)t$ where $x$ runs over $F_S$ and $t$ runs over
the maximal torus $T(F_S)$ of $GL_2(F_S)$ consisting of diagonal matrices. We 
substitute the Fourier expansion for $\varphi$ into (\ref{step1}).  We then move the $Y$-unipotents
and the matrix in $x$ leftward.
Since
$$w_{32} u_1(x) u_3(x)=u_{123}(x) u_2(x) w_{32}$$
and $u_{123}(x)$ is in the unipotent radical of $P$, the resulting integrand becomes
\begin{multline*}\sum_m\sum_{j,l} X_{lj}^{[d_1^{-1}d_2^{-1}]}
\frac{b_l(m)\,(m,d_1^{-1}d_2^{-1})_S}{|m|^{1/2}}\, \psi\left(\frac{mb_2c_1}{d_2}\right)
W_j\left(\left(\begin{smallmatrix}m&\\&1\end{smallmatrix}\right)
t\right)\\
f_s\left(\begin{pmatrix}d_1^{-1}d_2^{-1}&&&\\&d_1^{-1}d_2^{-2}&&\\&&d_1^{-2}d_2^{-1}&\\&&&1\end{pmatrix}
\begin{pmatrix}1&y_1&y_2&\\&1&x&\\&&1&\\&&&1\end{pmatrix}
w_{32}\begin{pmatrix}1&&&\\&1&x_1&x_2\\&&1&\\&&&1\end{pmatrix}
\begin{pmatrix}t&\\&t\end{pmatrix}\right)\\
\psi\left( mx+h(y_1+x_2) \right).
\end{multline*}

We now recognize the integral over $y_1,y_2$ and $x$ as being a Whittaker integral of the inducing data.
Denote the Whittaker coefficients of the theta function on the $n$-fold cover of $GL_3$ by $\tau(r_1,r_2)$,
normalized as in \cite{B-H2}, Section 4 (see that reference for more information about the coefficients).
Then we have expressed $I(s,\varphi,f_s)$ as a finite sum
\begin{equation}\label{finite-sum}
I(s,\varphi,f_s)=\sum_j D_j(s,\varphi)\, I_j(s,\varphi,f_s).
\end{equation}
Here $D_j(s,\varphi)$ is the Dirichlet series
\begin{multline*}
D_j(s,\varphi)=|h|^{-1}\!\!\!\sum_{0\neq m\in\mathcal{O}_S/\mathcal{O}_S^\times}\!\!
\sum_{\substack{0\neq d_1,d_2\in\mathcal{O}_S/\mathcal{O}_S^\times\\d_2|h,~d_1| md_2\\
c_1 \text{~mod~} d_1d_2,~(c_1,d_1)=1\\ c_2\text{~mod}^\times d_2}}\sum_l
\left(\frac{c_1}{d_1}\right)\left(\frac{c_2}{d_2}\right)\psi\left(\frac{h c_1}{d_1 d_2}\right)\,
 \psi\left(\frac{m b_2 c_1}{d_2}\right) \\ X_{lj}^{[d_1^{-1}d_2^{-1}]}\,b_l(m)\,(d_1d_2,m)_S\,
\tau(md_1^{-1}d_2,hd_2^{-1})
|d_1|\,|m|^{-2s+1/2} |d_1d_2|^{-4s}.
\end{multline*}
The divisibility conditions on $d_1$ and $d_2$ follow from the invariance of $f_s$ under lower
triangular matrices in $P(\mathcal{O}_S)$, and we have used that $(m,d_1^{-1}d_2^{-1})_S=(d_1d_2,m)_S$.
The factor that multiplies the series $D_j(s,\varphi)$ is 
\begin{multline*}I_j(s,\varphi,f_s)=\int_{X(F_S)}\int_{Z^n(F_S)\backslash T(F_S)}
W_j\left(t\right)\\
W_{f_s}\left(w_{32}\begin{pmatrix}1&&&\\&1&x_1&x_2\\&&1&\\&&&1\end{pmatrix}
\begin{pmatrix}t&\\&t\end{pmatrix}\right)
\psi\left(h x_2 \right)\,\delta_B(t)^{-1}\,dt\,dx_1\,dx_2,
\end{multline*}
where 
$$W_{f_s}(g)=\int_{N_3(\mathcal{O}_S)\backslash N_3(F_S)} f_s\left(
\begin{pmatrix} n&0\\0&1\end{pmatrix}g\right)\,\psi(n_{12}+n_{23})\,dn,$$
with $N_3$ the unipotent radical of standard Borel in $GL_3$.
This function may be expressed in terms of the Whittaker function at places in $S$ for $\Theta$.

Since the Eisenstein series has analytic continuation, we conclude that

\begin{theorem}  The function $\sum_j D_j(s,\varphi)I_j(s,\varphi,f_s)$ has analytic continuation to all complex $s$.
\end{theorem}

Of course, since the Eisenstein series has functional equation,
the series function $\sum_j D_j I_j$ does too.  But we shall not compute the
intertwining operators necessary to give this explicitly.

When $n=2$, $\Theta$ does not have a Whittaker model.
Since all non-degenerate Whittaker
coefficients $\tau(n_1,n_2)$ are zero, equation (\ref{finite-sum}) shows that the integral is zero,
confirming Conjecture~\ref{conje1}, Part 1, in this case.  (When $n=1$ there is no theta function, but
if one used the constant function in its place once again the integral would be zero by the same argument.)
We consider the first nonzero cases $n=3,4$ in the next Section.

\section{Analysis of $D_j(s,\varphi)$ when $n=3$ or $n=4$}\label{analysis}

As mentioned above,
Kazhdan and Patterson \cite{K-P} showed that the theta representation $\Theta$ on an $n$-fold cover of $GL_3(F_S)$
has a unique Whittaker model only when $n=3$ or $n=4$.  (If $n=4$ this
is only true for certain covers.)   It is a difficult problem to determine the Whittaker coefficients $\tau(r_1,r_2)$
of functions in the space of $\Theta$
for $n>4$; even on $GL_2$ the analogous coefficients are not determined (there for $n>3$).
In this Section we consider the two cases where one does have a unique Whittaker model and
use this information to analyze the series $D(s,\varphi)$.   

For convenience we suppose that $h$, which appears in the character $\psi_h$, equals $1$ (working with a more
general $h$ only changes the series at the primes dividing $h$).  Since $d_2|h$, we have $d_2=1$ as well.  The series
$D_j(s,\varphi)$ reduces to
\begin{multline*}
D_j(s,\varphi)=\\ \sum_{0\neq m\in\mathcal{O}_S/\mathcal{O}_S^\times} 
\sum_{\substack{0\neq d_1\in\mathcal{O}_S/\mathcal{O}_S^\times\\
d_1| m}} \sum_l g(d_1)\,(d_1,m)_S\,X_{lj}^{[d_1^{-1}]}\,
b_l(m)\,  \tau(md_1^{-1},1) \,
|m|^{-2s+1/2} |d_1|^{1-4s},
\end{multline*}
where $g(d_1)$ is the Gauss sum (for $\mathcal{O}_S$)
$$g_t(d_1)=\sum_{c_1 \text{~mod}^\times d_1}\left(\frac{c_1}{d_1}\right)^t\,\psi\left(\frac{c_1}{d_1}\right)$$
with $t=1$.
After an interchange of summation (and using properties of the Hilbert symbol and also replacing $d_1$ by simply $d$), we obtain
\begin{multline}\label{niceseries}
D_j(s,\varphi)= \\ 
\sum_{0\neq d\in\mathcal{O}_S/\mathcal{O}_S^\times}
\sum_{0\neq m\in\mathcal{O}_S/\mathcal{O}_S^\times}\sum_l X_{lj}^{[d^{-1}]}
g(d)\,(m,d)_S\,
b_l(md)\,  \tau(m,1) \,|m|^{-2s+1/2} |d|^{3/2-6s}.
\end{multline}

We recall that the Gauss sum $g(d)$ satisfies the equation
$$g(d_1d_2)=\left(\frac{d_1}{d_2}\right)\left(\frac{d_2}{d_1}\right)\,g(d_1)\,g(d_2)$$
provided $(d_1,d_2)=1$, and that $g(p^j)=0$ for any prime $p$ if $j\geq2$.

Consider first the case $n=3$.  The coefficients $\tau(r_1,r_2)$ were computed by Proskurin \cite{Pr}
and by Bump and Hoffstein \cite{B-H0}.  In particular, $\tau(m,1)=0$ unless $m$ is a unit times a perfect cube,
and $\tau(a m^3,1)=|m|\,\tau(a,1)$ for $0\neq a\in\mathcal{O}_S$.  After
some relabeling of variables, the series
then becomes
\begin{equation*}
D_j(s,\varphi)=  
\sum_{0\neq d\in\mathcal{O}_S/\mathcal{O}_S^\times}
\sum_{0\neq m\in\mathcal{O}_S/\mathcal{O}_S^\times}\sum_l X_{lj}^{[d^{-1}]}
g(d)\,b_l(m^3d)\,  \,|m|^{-s_1} |d|^{-s_1-1},
\end{equation*}
where $s_1=6s-5/2$.

Suppose that $\varphi$ is a Hecke eigenform at $p$ with eigenvalue $\lambda_p$.  Then if $(p,M)=1$,
the Fourier coefficients $b_l(m)$ satisfy a relation
(\cite{B-H2}, Cor.\ 3.3)
\begin{multline*}\lambda_p b_l(p^k M)= \\
b_l( p^{k-3}M)+b_l(p^{k+3}M) + |p|^{-1}\left(\frac{p}{M}\right)^{-k-1}g_{k+1}(p)\,
 \sum_{i} X_{il}^{[p^{-k-1}]} b_i(p^{1-k}M).
 \end{multline*}
Here the coefficients $b_l(m)$ are zero by definition if $m$ is not in $\mathcal{O}_S$.

To analyze the contribution obtained from $D_j(s,\varphi)$ at $p$, we compute
$$(1-\lambda_p p^{-s_1} + p^{-2s_1})\,D_j(s,\varphi).$$
Since $g(d)=0$ if $p^2|d$, there are two nonzero contributions, from $d$ such that $(d,p)=1$
and from $d$ such that $\text{ord}_p(d)=1$.  For the first contribution, the Hecke relation can
be rewritten
$$\lambda_p b_l(p^{3k}M)=\begin{cases}
b_l(p^{3(k-1)}M)+b_l(p^{3(k+1)}M)&\text{if $k\geq1$}\\
b_l(p^3M)+|p|^{-1}\left(\frac{p}{M}\right)^{-1}{g}(p)\sum_{i} X_{il}^{[p^{-1}]} b_i(pM)&
\text{if $k=0$}\end{cases}$$
where $(p,M)=1$.
Applying this it is not difficult to see that
\begin{multline}\label{first-term}
(1-\lambda_p p^{-s_1} + p^{-2s_1})\times \\
\sum_{\substack{0\neq d\in\mathcal{O}_S/\mathcal{O}_S^\times\\ (d,p)=1}}
\sum_{0\neq m\in\mathcal{O}_S/\mathcal{O}_S^\times}\sum_l X_{lj}^{[d^{-1}]}
g(d)\,b_l(m^3d)\,  \,|m|^{-s_1} |d|^{-s_1-1}
=\\
\sum_{\substack{0\neq d\in\mathcal{O}_S/\mathcal{O}_S^\times\\ (d,p)=1}}
\sum_{\substack{0\neq m\in\mathcal{O}_S/\mathcal{O}_S^\times\\ (m,p)=1}}\sum_l X_{lj}^{[d^{-1}]}
g(d)\,|m|^{-s_1} |d|^{-s_1-1}\\
\left(b_l(m^3d) - \left(\frac{p}{d}\right)^{-1}g(p)\sum_i X_{il}^{[p^{-1}]}b_i(pm^3d)|p|^{-s_1-1}\right).
\end{multline}
To analyze the terms with $\text{ord}_p(d)=p$, we make use of the Hecke relation 
$$\lambda_p b_l(p^{3k+1}M)=\begin{cases}
\bar\chi(p)\,b_l(p^{3(k-1)+1}M)+b_l(p^{3(k+1)+1}M)&\text{if $k\geq1$}\\
b_l(p^4M)+|p|^{-1}\left(\frac{p}{M}\right)\bar{g}(p)\sum_{i} X_{il}^{[p^{-2}]} b_i(M)&
\text{if $k=0$}\end{cases}$$
where $(p,M)=1$.  We also recall the properties (\cite{B-H2}, Eqs.\ (3.11) and (3.12))
\begin{equation}\label{x-behavior1}
X_{jk}^{[ab]}=(b,a)_S\,\sum_\ell X_{j \ell}^{[a]}X_{\ell k}^{[b]}\qquad a,b\in F^\times
\end{equation}
and 
\begin{equation}\label{x-behavior2}
X_{jk}^{[p^3]}=\delta_{j,k}\qquad \text{Kronecker delta.}
\end{equation}

Using the Hecke relations, Eqns.\ (\ref{x-behavior1}) (with $a=p^{-1}$, $b=d^{-1}$)
and (\ref{x-behavior2}) and reordering, one finds that 
\begin{multline}\label{second-term}
(1-\lambda_p p^{-s_1} + p^{-2s_1}) \times\\
\sum_{\substack{0\neq d\in\mathcal{O}_S/\mathcal{O}_S^\times\\ (d,p)=1}}
\sum_{0\neq m\in\mathcal{O}_S/\mathcal{O}_S^\times}\sum_l X_{lj}^{[d^{-1}p^{-1}]}
g(dp)\,b_j(m^3dp)\,  \,|m|^{-s_1} |dp|^{-s_1-1}
=\\
\sum_{\substack{0\neq d\in\mathcal{O}_S/\mathcal{O}_S^\times\\ (d,p)=1}}
\sum_{\substack{0\neq m\in\mathcal{O}_S/\mathcal{O}_S^\times\\ (m,p)=1}}
\sum_\ell X_{\ell j}^{[d^{-1}]}
g(d)\,g(p)\left(\frac{p}{d}\right)\left(\frac{d}{p}\right)\,|m|^{-s_1} |d|^{-s_1-1}\\
\left((d,p)_S\,\sum_l X_{l\ell}^{[p^{-1}]} b_l(pm^3d)\,|p|^{-s_1-1}-|p|^{-2s_1-2}\left(\frac{p}{d}\right)\bar{g}(p)
b_\ell(md^3)\right).
\end{multline}
We may now combine the terms (\ref{first-term}) and (\ref{second-term}).  Indeed, after using cubic
reciprocity, we see that the last term on the right-hand side
of (\ref{first-term}) cancels a term on the right-hand-side of (\ref{second-term}).  
Then iterating over the primes $p$ of $\mathcal{O}_S$, we obtain
the following Theorem.

\begin{theorem}\label{euler-product}  Suppose that $n=3$.  Let $$D_j^*(s,\varphi)=\zeta(12s-4)\,D_j(s,\varphi).$$
Then for each $j$, 
$$D^*_j(s,\varphi)=b_j(1)\,\prod_{p}\left(1-\lambda_p |p|^{-s_1} + |p|^{-2s_1}\right)^{-1}$$
with $s_1=6s-5/2$.
\end{theorem}
Here for arbitrary $n\geq2$, $\zeta(4ns-2n+2)$ is the normalizing factor of the Eisenstein series
$E_\Theta(g,s,f_s)$.

When $n=3$, the series considered by Bump and Hoffstein \cite{B-H2} is essentially of the form
$$\sum_l\sum_{m_1,m_2} \bar{\tau}(m_1,m_2)\,b_l(m_1) |\,m_1m_2^2|^{-s}.$$
Though it represents the same Euler product (at $3s-1$), the convolution (\ref{niceseries}) given
here is visibly different.

To conclude this section, we briefly discuss the case $n=4$.  In this case, $\varphi$ is an automorphic
form on the {\it double} cover of $GL_2(F_S)$. In (\ref{niceseries}), the only coefficients  $\tau(m,1)$
of the theta function on the 4-fold cover of $GL_3(F_S)$ that are not zero
 are those that occur when each prime dividing $m$ occurs to a power congruent to $0$
or 1 modulo $4$.   Moreover,  at prime powers we have
$$\tau(p^{4k},1)=|p|^k\qquad \tau(p^{4k+1},1)=|p|^{k-1/2} \bar{g}(p).$$
Thus the series involves both the Fourier coefficients of $\varphi$ at primes $p$, whose squares are related to the central
values of the Shimura lift of $\varphi$ twisted by a quadratic character modulo $p$, and quartic
Gauss sums $g(p)$.

\bigskip
\noindent
\textsc{\small{Department of Mathematics, Boston College, Chestnut Hill, MA
02467-3806, USA}}

\medskip
\noindent
\textsc{\small{School of Mathematical Sciences, Tel Aviv University, Ramat Aviv, Tel Aviv 6997801,
Israel}}

\end{document}